\documentclass{amsart}

\usepackage{hyperref}
\usepackage{url}

\usepackage{enumerate}        %for easier lists

\newtheorem{theoremIntro}{Theorem}

\renewcommand{\int}{\operatorname{int}}
\newtheorem{theorem}{Theorem}[section]
\newtheorem{lemma}[theorem]{Lemma}
\newtheorem{question}[theorem]{Question}

\newtheorem{corollary}[theorem]{Corollary}

\theoremstyle{definition}
\newtheorem{definition}[theorem]{Definition}

\newtheorem{remark}[theorem]{Remark}

\newcommand{\Z}{\mathbb{Z}}
\def\pol{p^X}
\def\qol{q^X}

\newcommand{\link}{\operatorname{link}}
\newcommand{\tu}{\mathbf{t}}
\newcommand{\turev}{\overline{\mathbf{t}}}

\begin{document}

\title{Face vectors of subdivided simplicial complexes}
\author{Emanuele Delucchi}
\author{Aaron~Pixton}
\author{Lucas~Sabalka}

\begin{abstract}

Brenti and Welker have shown that for any simplicial complex $X$, the face vectors of successive barycentric subdivisions of $X$ have roots which converge to fixed values depending only on the dimension of $X$.  We improve and generalize this result here.  We begin with an alternative proof based on geometric intuition.  We then prove an interesting symmetry of these roots about the real number $-2$.  This symmetry can be seen via a nice algebraic realization of barycentric subdivision as a simple map on formal power series in two variables.  Finally, we use this algebraic machinery with some geometric motivation to generalize the combinatorial statements to arbitrary subdivision methods:  any subdivision method will exhibit similar limit behavior and symmetry.  Our techniques allow us to compute explicit formulas for the values of the limit roots in the case of barycentric subdivision.
\end{abstract}

\maketitle

%%%%%%%%%%%%%%%%%%%%%%%%%%%%%%%%%%%%%%%%%%%%%%%%%%%
\section{Introduction}\label{sec:intro}

 Throughout this paper, we let $X$ be an arbitrary finite simplicial complex of dimension $d-1$, and we assume that all vectors and matrices will be indexed by rows and columns starting at $0$.  We are interested in roots of the \emph{$f$-polynomial} of $X$, defined as follows.  Let $f^X_i$ denote the number of $i$-dimensional faces of $X$.  We declare that $f^X_{-1} = 1$, where the $(-1)$-dimensional face is the empty face, $\emptyset$.  The \emph{face vector}, or \emph{$f$-vector}, of $X$ is the vector 
	$$f^X := (f^X_{-1},  f^X_0, \dots, f^X_{d-1}).$$
Let $\tu$ denote the column vector of powers of $t$, $(t^d,t^{d-1},\dots t^0)^T$.  The \emph{$f$-polynomial} $f^X(t)$ encodes the $f$-vector as a polynomial:
	$$f^X(t) := \sum_{j=0}^d f^X_{j-1} t^{d-j} = f^X\tu.$$

Much work has been devoted to the study of  $f$-vectors of simplicial complexes, their close relatives, the $g$- and $h$-vectors, and the associated polynomials.  As it turns out, the entries of these objects encode many combinatorial and algebraic aspects of the complex to which they are associated (see \cite{BayerLee, KleeKlein, Ziegler} for background and further references).

We focus on a recent result of Brenti and Welker which may initially appear surprising.  Let $X'$ denote the barycentric subdivision of $X$, and more generally let $X^{(n)}$ denote the $n^{\textrm{th}}$ barycentric subdivision of $X$. 

\begin{theorem}\cite{BrentiWelker}\label{thm:BrentiWelker}
Let $X$ be a $(d-1)$-dimensional simplicial complex.  As $n$ grows, the $d-1$ largest roots of $f^{X^{(n)}}(t)$ converge to $d-1$ negative real numbers which depend only on $d$, not on $X$.
\end{theorem}

We provide some geometric intuition and motivation for why this result holds.  We offer an alternate proof of this theorem based on these geometric observations.  In the process, we show how to compute the $d-1$ real values for each $d$.  Our first main theorem is:

\begin{theoremIntro}\label{thm:main}
Let $X$ be a $(d-1)$-dimensional simplicial complex.  Then the $d-1$ largest roots of $f^{X^{(n)}}(t)$ converge to $d-1$ values which are the roots of a polynomial $p_d(t)$, depending only on $d$, whose coefficients are listed in the last row of the inverse of a particular matrix, $P_d$.
\end{theoremIntro}
The entries of the matrix $P_d$, and of its inverse $P_d^{-1}$, are computed in Section \ref{sec:computations}.   Our calculations allow us to compute the `limit roots' thus obtained.  In the examples, we observed that these `limit roots' are symmetrically distributed about the point $-2$, with respect to the M\"{o}bius transformation $x \mapsto \frac{-x}{x+1}$.  Our second main theorem proves this symmetry:

\begin{theoremIntro}\label{thm:symmetry}
For any dimension $(d-1)$, the $d-1$ `limit roots' are invariant under the map $x \mapsto \frac{-x}{x+1}$.
\end{theoremIntro}

In fact, more can be said.  The existence of a `limit polynomial' and the symmetry result hold for an arbitrary subdivision method, as we show in Theorem \ref{thm:generalsym} .\\
%(see Section \ref{sec:geometry} for definitions):

%\begin{theorem}\label{thm:generalsymmetry} (see Theorem \ref{thm:generalsym})
%For any dimension $n$ and any subdivision method $\Phi$ which is nontrivial in dimension $n$, there exists a unique `limit polynomial' $p_{n,\Phi}(t)$.  The polynomial $p_{n,\Phi}(t)$ is such that, for any $d$-dimensional simplicial complex $X$, the roots of the $f$-polynomials obtained by repeatedly subdividing $X$ converge to the roots of $p_{n,\Phi}(t)$.  The roots of $p_{n,\Phi}(t)$ are invariant under the M\"{o}bius transformation $x \mapsto \frac{-x}{x+1}$.
%\end{theorem}

For barycentric subdivision, this symmetry can be seen through a beautiful algebraic theorem.  Barycentric subdivision, considered as a map on $f$-polynomials, induces a function $b: \Z[t] \rightarrow \Z[t]$, as in Section \ref{sec:sym}.  We list the values of $b$ on monomials as coefficients in the formal power series in the variable $x$ over $\Z[t]$, by defining $B:\Z[t][[x]] \to \Z[t][[x]]$ by $B(\sum_{k\geq 0}g_k(t) x^k) := \sum_{k\geq 0}b\big(g_k(t)\big)x^k$.

\begin{theoremIntro}\label{thm:powerseries}
In $\mathbb{Z} [t][[x]]$, barycentric subdivision satisfies the identity
$$B(e^{tx})=\frac{1}{1-(e^x-1)t}.$$
\end{theoremIntro}

This paper is organized as follows.  In Section \ref{sec:Motivation}, we discuss the geometric intuition and motivation behind Theorem \ref{thm:BrentiWelker}.  In Section \ref{sec:main}, we prove Theorem \ref{thm:main}.  In Section \ref{sec:sym}, we prove the symmetry stated in Theorem \ref{thm:symmetry}, and prove Theorem \ref{thm:powerseries}.  In Section \ref{sec:generalsub}, we extend the symmetry to arbitrary subdivision methods.  We end with Section \ref{sec:computations}, where we compute the entries in the inverse matrix $P_d^{-1}$ found in Theorem \ref{thm:main} as well as all limit polynomials and roots up to the value $d = 10$.

\subsection*{Acknowledgements.} The authors would like to thank Laura Anderson for bringing Theorem \ref{thm:BrentiWelker} to their attention, Thomas Zaslavsky for interesting discussions, and Dennis Pixton for introducing the authors to each other.

%%%%%%%%%%%%%%%%%%%%%%%%%%%%%%%%%%%%%%%%%%%%%%%%%%%
\section{Geometric Motivation}\label{sec:Motivation}

Brenti and Welker's theorem may be surprising at first:  there is no dependence on the initial complex $X$, only on the dimension $d-1$.  However, geometrically this makes perfect sense.  Barycentrically subdividing a simplicial complex $X$ over and over again causes the resulting complex $X^{(n)}$ to have far more cells than the original $X$.  Because higher-dimensional cells contribute more new cells (in every dimension) upon subdividing than lower-dimensional ones, the top-dimensional cells begin to dominate in their number of contributions to subdivisions.  For example, think of geometric realizations so that $X^{(n)}$ is a subset of $X$.  Then a randomly chosen cell of $X^{(n)}$ should, with higher and higher probability as $n$ increases, be contained in the interior of a top-dimensional cell of $X$, as top-dimensional cells contribute far more cells to $X^{(n)}$ than other cells.  

Each of the $f^X_{d-1}$ top-dimensional cells of $X$ contributes the same number of cells to $X^{(n)}$.  Since these cells eventually dominate contributions from smaller-dimensional cells, the $f$-polynomial for $X^{(n)}$ can be approximated by $f^X_{d-1}$ times the $f$-polynomial for the $n^{th}$ barycentric subdivision $\sigma_d^{(n)}$ of a single top-dimensional cell $\sigma_d$.  Since the roots of a polynomial are unaffected by multiplication by constants, the roots of $f^{X^{(n)}}(t)$ converge to the roots of $f^{\sigma_d^{(n)}}(t)$ as $n$ increases.  

By definition, the coefficients of $f^{\sigma_d^{(n)}}(t)$ records the number of cells of each dimension occurring in $\sigma_d^{(n)}$.  The number of cells in each dimension is bounded by a constant times the number of top-dimensional cells.  Thus, if we normalize $f^{\sigma_d^{(n)}}(t)$ by dividing by the number of top-dimensional cells, we have coefficients which, for each $k$, record the density of $k$-cells relative to the number of top-dimensional cells.  As this density is positive but strictly decreases upon subdividing, there is a limiting value for the coefficient.  Thus, there is a limiting polynomial, with well-defined roots. 

We now formalize this intuition.

%%%%%%%%%%%%%%%%%%%%%%%%%%%%%%%%%%%%%%%%%%%%%%%%%%%
\section{$f$-polynomials of barycentric subdivisions} \label{sec:main}

%%%%%%%%%%%%%%%%%%%%%%
\subsection{Barycentric Subdivision and the Matrix $\Lambda_d$}
To prove Theorem \ref{thm:main}, we begin by observing the effect of barycentric subdivision on $f$-vectors.  One key observation is that barycentric subdivision multiplies $f$-vectors by a fixed matrix, $\Lambda_d$, defined as follows.

\begin{definition} For $i,j \geq -1$, let $\lambda_{i,j}$ denote the number of $j$-dimensional faces in the interior of the first barycentric subdivision of the standard $i$-dimensional simplex, where by convention $\lambda_{-1,-1} = 1$ and $\lambda_{i,-1}$ is $1$ if $i = -1$ and $0$ otherwise.  Let $\Lambda_d$ denote the $(d+1) \times (d+1)$ matrix, with rows and columns indexed by the integers $-1, 0, \dots, d-1$, whose entry in the $i^{th}$ row and $j^{th}$ column is $\lambda_{i,j}$:  $\Lambda_{d} := [\lambda_{i,j}]$.
\end{definition}

For example, $\Lambda_5 = \left[\begin{array}{cccccc} 
1 & 0 & 0 & 0 & 0 & 0\\
0 & 1 & 0 & 0 & 0 & 0\\
0 & 1 & 2 & 0 & 0 & 0\\
0 & 1 & 6 & 6 & 0 & 0\\
0 & 1 & 14 & 36 & 24 & 0\\
0 & 1 & 30 & 150 & 240 & 120
\end{array}\right]$

\begin{lemma} Barycentrically subdividing a $(d-1)$-dimensional simplicial complex $X$ multiplies the $f$-vector by $\Lambda_d$:
	$$
	f^{X'}= f^X\Lambda_d.
	$$
\end{lemma}

\begin{proof} The faces of $X'$ can be partitioned according to the lowest-dimensional faces of $X$ containing them.  Each face of $X$ is a simplex of some dimension $i$, and thus its interior contributes $\lambda_{i,j}$ to the total number of $j$-cells of $X'$ (or, if $i = 0$, exactly one vertex to $X'$). The claim then follows by linearity.
\end{proof}

\begin{corollary}\label{cor:fXn}
For any $n \geq 0$, 
	$$f^{X^{(n)}}=f^X \Lambda_d^n.$$
\end{corollary}

Thus, to understand barycentric subdivision, we need to understand the matrix $F_d$.  We will compute the entries in $F_d$ more explicitly in the following two sections, but for now we simply observe a formula, also noted in \cite{BrentiWelker}, which follows from Inclusion-Exclusion:

\begin{lemma}\label{lem:foij}
If $j > i$ then $\lambda_{i,j}= 0$.  If $j \leq i$, then
	$$\lambda_{i,j}=\sum_{k=0}^{i+1} (-1)^k{i+1 \choose k } f_j^{\sigma'_{i+1-k}}. \qquad\qed $$ 
\end{lemma}

By this lemma, $F_d$ is lower triangular with diagonal entries $\lambda_{i,i} =f_i^{\sigma'_i}=i!$.  Thus, the eigenvalues of $F_d$ are $0!, 1!, 2!, 3!, \dots, d!$.

%%%%%%%%%%%%%%%%%%%%%%
\subsection{Limit Behavior of the Roots} We now turn to the roots of the $f$-polynomials $f^{X^{(n)}}(t)$.  

Note that some results in this section were exhibited in or follow from \cite{BrentiWelker}.  We present here a self-contained exposition in order to present an alternative and more geometrically motivated proof of (most of) Brenti and Welker's result, Theorem \ref{thm:BrentiWelker}.

By Corollary \ref{cor:fXn}, 
	$$
	f^{X^{(n)}}(t)= f^X \Lambda_d^n \tu .
	$$
As the greatest eigenvalue of $\Lambda_d$ is $d!$, we normalize $f^{X^{(n)}}(t)$ by dividing by $(d!)^n$ - let $\pol_n(t)$ denote the result:
	$$
	\pol_{n}(t) := \frac{1}{(d!)^n}f^{X^{(n)}}(t).
	$$
Note that this normalization does not alter the roots.  It will also often be convenient to reverse the order of the coefficients of $\pol_n(t)$, with the effect of inverting the roots of $\pol_n(t)$ (that is, the roots of $f^{X^{(n)}}(t)$) about the unit circle in the extended complex plane:
	$$
	\qol_{n}(t) := t^d\pol_{n}(t^{-1}).
	$$

We are interested in the behavior of the roots of $\pol_n(t)$ and $\qol_n(t)$ as $n$ goes to infinity, so we are interested in the powers of $\Lambda_d$.  To take powers of $\Lambda_d$, we diagonalize,
	$$\Lambda_d=P_dD_dP_d^{-1},$$
where $D_d$ is the diagonal matrix of eigenvalues $0!, 1!, \dots, d!$ and $P_d$ is the (lower triangular) diagonalizing matrix of eigenvectors.  Thus, $\Lambda_d^n = P_dD_d^nP_d^{-1}$.  

Let $\widetilde{D}_d:=\frac{1}{d!}D_d$.   Let $\turev$ denote the column vector $\tu$ in reverse order, $\turev = (t^0,t^1,\dots t^d)^T$.  For any simplicial complex $X$, we thus have the following equations:
	\begin{eqnarray*}
	f^{X^{(n)}}(t) 
	&=& f^X P_d D_d^n P_d^{-1} \tu \\
	&=& (d!)^n \left(f^X  P_d\right) \left(\widetilde{D}_d\right)^n \left(P_d^{-1} \right) \tu,\\\\
	\pol_{n}(t) 
	&=& \left(f^X  P_d\right) \left(\widetilde{D}_d\right)^n \left(P_d^{-1} \right) \tu,\\
	\qol_{n}(t) 
	&=& \left(f^X  P_d\right) \left(\widetilde{D}_d\right)^n \left(P_d^{-1} \right) \turev.
	\end{eqnarray*}

The goal of Section \ref{sec:computations} will be to describe more precisely the matrices $P_d$ and $P_d^{-1}$.  As the eigenvalues of $\Lambda_d$ are $0!, 1!, \dots, d!$, for large $n$, $D_d^n$ is dominated by its last diagonal entry, $(d!)^n$.  In the limit, the powers of the matrix $\widetilde{D}_d=\frac{1}{d!}D_d$ converge to the matrix
	$$M_{d,d} := 
	\begin{bmatrix}
		0		& \cdots	& 0		\\
		\vdots	& \ddots	& \vdots	\\ 
		0		& \cdots	& 1
	\end{bmatrix}.$$
Thus, as $n$ grows, the polynomials $\pol_{n}$ and $\qol_n$ approach the polynomials
	\begin{eqnarray*}
	\pol_{\infty}(t) &:=& \left(f^X  P_d\right) M_{d,d} \left(P_d^{-1} \right) \tu \text{ and}\\
	\qol_{\infty}(t) &:=& \left(f^X  P_d\right) M_{d,d} \left(P_d^{-1} \right) \turev,
	\end{eqnarray*}
respectively, in the sense that each sequence converges coefficient-wise in the vector space of polynomials of degree at most $d$. 

By Corollary \ref{cor:fXn} and Lemma \ref{lem:foij}, we know the leading and trailing coefficients of $\pol_n(t)$ and $\qol_n(t)$:  $\pol_n(t)= (d!)^{-n}t^{d} + \dots + f^X_{d-1}$ and $\qol_n(t) = (d!)^{-n}  + \dots + f^X_{d-1}t^{d}$.  Hence, in the limit, $\pol_\infty(t)$ does not have $0$ as a root, but has degree less than $d$ (one root of the $\pol_n$ diverges to $-\infty$), while $\qol_\infty(t)$ is of degree $d$ with $0$ as a root.  Because the polynomials $\qol_n(t)$ converge coefficient-wise to the polynomial $\qol_\infty(t)$ of the same degree, their roots also converge:

\begin{lemma}\cite{BookMasur}\label{lem:rootsconverge}
 Let $(P_n(t))_n$ be a sequence of monic polynomials of degree $d$ that converges to a monic polynomial $P_\infty(t)$ of the same degree $d$. Then the roots of $P_n(t)$ may be numbered as $r^n_1,\dots,r_d^n$ and the roots of $P_\infty(t)$ as $r^\infty_1,\dots , r^\infty_d$ in such a way that for all $j=1,\dots , d$ the sequence $r_j^n$ converges to $r_j^\infty$ for $n\rightarrow \infty$. 
\end{lemma}

%\begin{proof}
%It is known that every root of a polynomial is a continuous function of the polynomial's coefficients (for a proof, see\cite{Tyrtyshnikov}). More precisely, given a polynomial $a_0+a_1t+\dots a_dt^d$, there are continuous functions $r_1(a_0,\dots ,a_d),\dots , r_d(a_0,\dots a_d)$ that parametrize the roots (with multiplicity). 
%\end{proof}

Since the roots of $\qol_n(t)$ converge to the roots of $\qol_\infty(t)$, it follows that the roots of $\pol_n(t)$ converge to the roots of $\pol_\infty(t)$ (with one of the roots `converging' to $-\infty$).

Because the matrix $P_d$ is lower triangular and $M_{d,d}$ has only one nonzero entry in position $(d-1,d-1)$, we have
	$$\left(f^X P_d\right) M_{d,d} = c_{X,d} e_d^T,$$
where $e_d$ is the $d^{th}$ unit vector, and $c_{X,d}$ is a constant depending on $f^X$ and $P_d$.  As both $f^X$ and $P_d$ do not depend on the amount of subdivision $n$, the roots of $\pol_\infty$ and $\qol_\infty$ do not depend on the value of $c_{X,d}$, and thus do not depend on \emph{any} coefficient of $f^X_d$.  This leads us to the following definition.

\begin{definition}
Define the \emph{limit $p$-polynomial} by
	$$p_d(t) := e_d^TP_d^{-1}\tu,$$
and the \emph{limit $q$-polynomial} by
	$$q_d(t) := e_d^TP_d^{-1}\turev.$$
\end{definition}

In this section we have proven:

\begin{theorem}\label{thm:facts}
The following facts hold:
\begin{enumerate}
\item The roots of $f^{X^{(n)}}(t)$ are equal to the roots of $\pol_n(t)$. 
\item The roots of $\qol_n(t)$ converge to the roots of $q_d(t)$, and depend only on the dimension of $X$.  
\item The roots of $\pol_n(t)$ converge to the roots of $p_d(t)$, and depend only on the dimension of $X$.  
\item The coefficient of $t^i$ in the polynomial $p_d(t)$ is the $(d-i)^{th}$ entry in last row of $P_{d+1}^{-1}$.  
\item The coefficient of $t^i$ in the polynomial $q_d(t)$ is the $(i-1)^{th}$ entry in the last row of $P_{d+1}^{-1}$.
\end{enumerate}
\end{theorem}

This proves Theorem \ref{thm:main}.  Note the first two facts give an alternative proof of Brenti and Welker's result, Theorem \ref{thm:BrentiWelker}, except for the fact that the roots are all real (that the roots are negative would then follow from the fact that all coefficients of these polynomials are positive).

In Section \ref{sec:computations}, we will explore the final two facts of Theorem \ref{thm:facts} by computing the entries of $P_d^{-1}$.

%%%%%%%%%%%%%%%%%%%%%%%%%%%%%%%%%%%%%%%%%%%%%%%%%%%
\section{Symmetry of the roots}\label{sec:sym}

Our goal is now to show that the limits of the roots satisfy the symmetry stated in Theorem \ref{thm:symmetry}.  We will prove this symmetry for the roots of $q_d$ instead of $p_d$, as it becomes a mirror symmetry instead of a M\"{o}bius invariance.

\begin{theorem}\label{thm:sym}
For every $d$, 
$$q_d(t)= (-1)^{d} q_d(-1-t).$$
In particular, the roots of $q_d(t)$ are (linearly) symmetric with respect to $-\frac{1}{2}$. 
\end{theorem}

To prove this theorem, we start by examining the subdivision of a single closed simplex.  The following lemma uses the usual difference operator $\Delta$ on a sequence, which takes a sequence $\{a_n\}_{n\geq 0}$ and returns the sequence $\{a_n-a_{n-1}\}_{n\geq 1}$.  We abuse notation by using $\Delta\{a_n\}_{n\geq 0}$ to also denote the first term in this sequence, with context determining whether the result is a single term or a sequence.

\begin{lemma}\label{lem:single}  Let $\sigma_s$ be a closed simplex of dimension $s-1$. The $f$-vector of the barycentric subdivision $\sigma_s'$ of $\sigma_s$ is given by
	$$f^{\sigma_s'}_j = \Delta^{d-j} \{f^{\sigma_s}(l)\}_l = \Delta^{d-j} \{(1+l)^s\}_l.$$
\end{lemma}

\begin{proof}
Lemma 2.1 in \cite{BrentiWelker} states that
\begin{align*}
f^{\sigma_s'}(t)
	&:=\sum_{j=0}^{s} f_{j-1}^{\sigma_s'}t^{d-j} \\
	&= \sum_{j=0}^{s} t^{d-j} \sum_{i=0}^{s} { s \choose i } \sum_{k=0}^{j}(-1)^k{j\choose k} (j-k)^i.\\
\end{align*}
Note that the innermost sum is the Stirling number $S(i,j)$ of the second kind (see \cite{Stanley}, page 34).  Reordering this triple summation, we have:
\begin{align*}
f^{\sigma_s'}(t)
	&= \sum_{j=0}^{s} t^{d-j} \sum_{k=0}^{j} (-1)^k{j\choose k} \sum_{i=0}^{s}  { s \choose i } (j-k)^i \\
	&= \sum_{j=0}^{s} t^{d-j} \sum_{k=0}^{j} (-1)^k{j\choose k} f^{\sigma_s}(j-k) \\
	&= \sum_{j=0}^{s} t^{d-j} \sum_{k=0}^{j} (-1)^{j-k}{j\choose k} f^{\sigma_s}(k) \\
	&=\sum_{j=0}^{s} \Delta^j \{f^{\sigma_s}(l)\}_l t^{d-j},
\end{align*}
where in the third equality we replace $k$ with $j-k$.
\end{proof}

\begin{corollary}
Let $X$ be a simplicial complex. The $f$-polynomial of its barycentric subdivision $f^{X'}(t)$ is given by 
	$$
	f^{X'}(t) = \sum_{j=0}^{d}\Delta^j \{f^X(l)\}_lt^{d-j}.
	$$
The polynomials $p_1^X$ and $q_1^X$ are given by
	$$
	(d!)\pol_1(t)=\sum_{j=0}^{d}\Delta^j \{\pol_0(l)\}_lt^{d-j} \quad \text{and} \qquad
	(d!)\qol_1(t)=\sum_{k=0}^{d}\Delta^k \{\qol_0(l)\}_lt^k. 
	$$
\end{corollary}
\begin{proof}
These formulas follow easily from Lemma \ref{lem:single} by linearity of the difference operator $\Delta$.
\end{proof}

Taking inspiration from the formula for $f^{X'}(t)$ above, we consider barycentric subdivision as a function on polynomials in $t$ defined by 
\begin{equation}\label{def:b}
b: \mathbb Z [t] \rightarrow \mathbb Z [t],\quad g(t) \mapsto \sum_{k\geq 0}\Delta^k \{g(l)\}_lt^k. 
\end{equation}
(Note that this sum is finite because the iterated finite differences of a polynomial are eventually all zero.)

For a simplicial complex $X$ of dimension $(d-1)$ we have that
$$
b(\qol_j(t))= d! \qol_{j+1}(t).
$$
The function $b$ is linear, and thus it is given by its values on monomials. It will be convenient to list these values as arranged on the `clothesline'~\cite{Wilf} provided by a formal power series in the variable $x$ over the ring $\mathbb Z [t]$. We thus consider a function $B$ on the ring $\mathbb Z [t] [[x]]$ defined as 
$$B:\quad \sum_{k\geq 0}g_k(t) x^k \quad \longmapsto \quad \sum_{k\geq 0}b\big(g_k(t)\big)x^k.$$

% It extends termwise to a function $B: \mathbb{Z} [[t]][[x]]\rightarrow \mathbb{Z} [[t]] [[x]]$  on formal power series over $\mathbb Z [[t]]$ with the understanding that the coefficient of $\frac{x^n}{n!}$ in the image of the formal power series for $e^{tx}$ will indicate the value of $b(t^n)$.

\begin{theorem} (see Theorem \ref{thm:powerseries})
In $ \mathbb{Z} [t] [[x]] $ it holds that
$$B(e^{tx})=\frac{1}{1-(e^x-1)t}.$$
\end{theorem}
\begin{proof}
We expand the right-hand side as a formal power series over $x$ and compare the coefficient of $\frac{x^n}{n!}$ therein with the value of $B(t^n)=b(t^n)$ as given in (\ref{def:b}). We have
\begin{align*}
\frac{1}{1-(e^x-1)t} &= \sum_{j\geq 0} (e^x-1)^jt^j = \sum_{j\geq 0} \bigg( \sum_{m=0}^j{ j \choose m}(-1)^{j-m}e^{mx} \bigg)t^j \\
&= \sum_{j\geq 0} t^j \bigg( \sum_{m=0}^j{ j \choose m}(-1)^{j-m}\sum_{k\geq 0}\frac{m^kx^k}{k!} \bigg)\\
&= \sum_{k\geq 0} \bigg(\sum_{j\geq 0} t^j \sum_{m=0}^j { j \choose m}(-1)^{j-m}m^k \bigg)\frac{x^k}{k!}\\
&= \sum_{k\geq 0} \bigg(\sum_{j\geq 0} \Delta^j\{m^k\}_{m} t^j \bigg)\frac{x^k}{k!}\\
&= \sum_{k\geq 0} b\big(t^k\big)\frac{x^k}{k!}= B\bigg(\sum_{k\geq 0} \frac{t^kx^k}{k!}\bigg) = B(e^{tx}).
\end{align*}
\end{proof}

To investigate the stated symmetry, we consider the following map
\begin{equation}\label{def:iota}
\iota:\mathbb Z [t] \rightarrow \mathbb Z [t], \quad g(t)\mapsto g(-1-t).
\end{equation}

\begin{lemma}\label{lem:invo}
The map $\iota$ is an involution, and it satisfies 
$$\iota b \iota  = b. $$
\end{lemma}
\begin{proof}
The map $\iota$ is clearly linear, so it will suffice to prove the claim for monomials.

It is easy to see that $\iota$ is an involution. Moreover, $\iota b = b \iota$, as
%\begin{align*}
\[
\iota B (e^{tx}) = \iota \bigg(\frac{1}{1-(e^{x}-1)t}\bigg) = \frac{1}{1-(e^{x}-1)(-t-1)}  = B(e^{(-1-t)x})  = B\iota (e^{tx}).
\]
%\end{align*}
The claim follows with term-by-term comparison.
\end{proof}

We are now ready to prove Theorem \ref{thm:sym}.
\begin{proof}[Proof of Theorem \ref{thm:sym}]
Barycentric subdivision has the effect on each $p$- and $q$-polynomial of multiplying on the right by $F$ before the $\tu$ and $\turev$, respectively, and rescaling by dividing by $d!$.  In the limit, the limit $p$- and $q$-polynomials are invariant under barycentric subdivision up to this scaling, so that
	$$
	b\big(q_d(t)\big)= d! q_d(t).
	$$
Since the eigenvalues of $F$ are all distinct, $q_d$ is characterized by this identity, and by having leading coefficient $f_{d-1}^X$.

Applying Lemma \ref{lem:invo}, we have
$$ 
 b \big(q_d(-1-t)\big)= b \big( \iota (q_d(t))\big) = \iota \big( b (q_d(t)) \big) =  \iota \big((d!) q_d(t)\big) = d! \big( q_d(-1-t) \big),
$$ 
and since the lead coefficient of $q_d(-1-t)$ is $(-1)^{d}f_{d-1}^X$, the claim follows.

\end{proof}

%%%%%%%%%%%%%%%%%%%%%%%%%%%%%%%%%%%%%%%%%%%%%%%%%%%
\section{Symmetry for Other Subdivision Methods} \label{sec:generalsub}

In general, given any polynomial $g(t) \in \Z[t]$, we can consider the polynomial $\iota g(t) = g(-1-t)$. The coefficient of $t^k$ in $g(t)$ contributes $(-1)^k{ k \choose j}$ times itself to the coefficient of $t^j$ in $\iota g(t)$: this contribution is up to sign the number of $(j-1)$-dimensional faces of the $(k-1)$-dimensional simplex.  Thus, we can interpret $\iota$ as a map on formal sums of simplices, as follows.

We will think of every simplex $\sigma\in X$ as a subset of the vertex set of $X$. Now we can write 
	$$
	\iota: \mathbb Z [X] \to \mathbb Z [X],\quad 
	\sigma\mapsto (-1)^{\dim \sigma+1}\sum_{\tau \subseteq \sigma} \tau.
	$$
We represent the simplicial complex $X$ as the formal sum $\sum_{\sigma\in X}\sigma$ of all its simplices, each with `weight' $1$.

Let us recall some basics about subdivisions of simplicial complexes, pointing to \cite{Spanier} as a reference for a more detailed discussion. In the following we will write $\vert X \vert$ for the geometric realization of a given simplicial complex $X$ \cite[Section 3.1]{Spanier}.

\begin{definition}[Compare Section 3.3 of \cite{Spanier}] A {\em subdivision} (not necessarily barycentric) of $X$ is a simplicial complex $\widetilde{X}$ whose vertices are points of $\vert X \vert$ and such that: 
\begin{enumerate}
\item For every simplex $\widetilde{\sigma}$ of $\widetilde{X}$ there is a simplex $\sigma$ of $X$ such that $\widetilde{\sigma}\subseteq \vert \sigma \vert$.
\item The linear map $\vert \widetilde{X} \vert \to \vert X \vert$ mapping each vertex of $\widetilde{X}$ to the corresponding point of $\vert X \vert$ is a homeomorphism.
\end{enumerate}

We will identify a subdivision of $X$ with the triple $(X,\widetilde{X}, \phi)$, where $\phi: \widetilde{X}\to X$ is the function associating to each $\widetilde{\sigma}$ the smallest simplex $\sigma \in X$ such that $\widetilde{\sigma} \subseteq \vert \sigma \vert$.

% triple $(X, \widetilde{X}, \phi)$ such that $\widetilde{X}$ is a simplicial complex with vertex set $V\widetilde{X}$, and $\phi:  V\widetilde{X} \to S$ is a map that, when extended to the set $\widetilde{S}$ of all simplices of $\widetilde {X}$ as $\phi: \widetilde{S}\to S,  \widetilde{\sigma} \mapsto\cup_{v \in \widetilde{\sigma}}\phi(v) $, satisfies the following conditions:
%	\begin{enumerate}
%	\item $S = \im(\phi)$  %$VX \subset \im(\phi)$.
%         \item For every simplex $\widetilde{\sigma} \subset V\widetilde{X}$ of $\widetilde{X}$ of dimension $\dim \widetilde{\sigma}$, the set $\phi(\widetilde{\sigma})$
%is a simplex of $X$ of dimension at least  $\dim \widetilde{\sigma}$.  %By abuse of notation, we write $\phi(\widetilde{\sigma}) =\cup_{v \in \widetilde{\sigma}}\phi(v)$.
%         \item There is a family of points $\{x(\widetilde{v})\}_{\widetilde{v}\in V\widetilde{X}}$ with  $x(\widetilde{v})\in \int (\vert \phi(\widetilde{v})\vert) \subseteq \vert X\vert$ such that the map $\widetilde{v}\mapsto x(\widetilde{v})$ extends linearly to a homeomorphism $\vert \widetilde{X}\vert \to \vert X\vert$.
%	\end{enumerate}
\end{definition}

Now, a subdivision $(X,\widetilde{X},\phi)$ induces a linear map
	$$
	b_{\phi}: \mathbb Z [X] \to \mathbb Z [\widetilde{X}],\quad 
	{\sigma} \mapsto \sum_{\phi(\widetilde{\sigma}) = \sigma}\widetilde{\sigma}.
	$$

In the following definition, we collect together compatible subdivisions in different dimensions, calling the result a \emph{subdivision method}.  This is not the most general definition of a subdivision method, and our results might hold in greater generality, but we restrict ourselves to the definition presented here because in more general subdivision methods can get notationally quite cumbersome without adding to the actual idea.

\begin{definition} A \emph{subdivision method} $\Phi$ is a collection of subdivisions $\Phi := \{(\sigma_n,\widetilde{\sigma}_n,\phi_n)\}_{n\geq 0}$ such that for every map $i_k:\sigma_k\to \sigma_m$ identifying a $k$-face of the standard $m$-simplex and every permutation of $m$ elements $\pi$, the map $\phi_k$ is the restriction of $\phi_m$ to $\pi i_k(\sigma_k)$.  This ensures that, given any simplicial complex $X$, the complex $\Phi (X)$, called {\em subdivision of $X$ according to the rule $\Phi$}, is uniquely defined by requiring that every $n$-simplex of $X$ is subdivided as $(\sigma_n,\widetilde{\sigma}_n,\phi_n)\in\Phi$.  A subdivision method is \emph{nontrivial in dimension $n$} if $\phi_k$ is not the identity map for some $k \leq n$.  Clearly if a subdivision is nontrivial in dimension $n$, then $\phi_n$ is not the identity map.
\end{definition}

Barycentric subdivision is the subdivision method where $\widetilde{\sigma}_n=2^{\sigma_n}$ and $\phi_n{j}=\{j\}$ for every vertex $j\in \sigma_n$.  

Given a subdivision method $\Phi$, in view of the linearity of $b_\phi$ for each subdivision, it makes sense to write 
	$$b_\Phi(\sum_{\sigma \in X} \sigma) = \sum_{\sigma \in X} b_{\Phi}\sigma.$$
As with the map $b$ given by barycentric subdivision, for any subdivision method the induced map $b_{\Phi}$ always commutes with the map $\iota$:

\begin{lemma}\label{lem:generalcommutativity}
For any subdivision method $\Phi$, $\iota b_\Phi = b_\Phi \iota$.
\end{lemma}

Before we prove this lemma, we need some properties of the map $\iota$.  For this paper, the \emph{link} of a simplex $\sigma$ in a simplicial complex $X$ is the subcomplex consisting of all simplices $\tau$ in $X$ such that $\sigma \cap \tau = \emptyset$ and $\sigma \cup \tau$, thought of as subsets of the vertex set of $X$, is also a simplex of $X$.  Applying $\iota$ to $X$, we obtain
	\begin{align}
	\notag\iota\bigg(\sum_{\sigma\in X}\sigma\bigg) &=
	\sum_{\sigma\in X} 
		(-1)^{\dim \sigma+1}\sum_{\tau\subseteq \sigma}\tau\\
         \notag&= \sum_{\tau\in X}\bigg(\sum_{\substack{\sigma\in X \\ \sigma\supseteq \tau}}(-1)^{\dim \sigma - \dim \tau}\bigg)(-1)^{\dim \tau+1}\tau\\
	 &= \sum_{\tau\in X} (-1)^{\dim \tau}\big(\chi(\link \tau )-1\big)\tau ,\label{eq:5}
	\end{align}
where $\chi$ is the Euler characteristic. That the Euler characteristic satisfies this identity can be found in \cite{StanleyCCA}.

We now need to characterize how $\iota$ acts on simplices.  We do so by looking at how $\iota$ affects a (not necessarily pure) simplicial homology manifold.  If $M$ is an $r$-dimensional (not necessarily pure) simplicial homology manifold with boundary $\partial M$ and $X$ is a finite simplicial complex such that $M$ is (PL-homeomorphic to) the geometric realization of $X$, then we let $[M]$ denote the formal sum of all simplices of $X$. In this setting the boundary submanifold $\partial M\subset M$ induces a subcomplex $\partial X\subset X$, and we take $[\partial M]$ to be the sum of all simplices in $\partial X$.

\begin{lemma}\label{lem:manifold formula}
For any (not necessarily pure) simplicial homology $r$-manifold $M$ with boundary $\partial M$ and dimension $r\ge 0$, 
	$$\iota([M]) = (-1)^{r+1}\big([M] - [\partial M]\big).$$
\end{lemma}

\begin{proof}
The link of every simplex $\sigma \in [M]$ is of dimension $r - \dim \sigma - 1$, and is a homology ball or sphere according to whether $\sigma$ is on the boundary $\partial M$ or not.  If $\link \sigma$ is a homology ball, $\chi (\link \sigma) - 1 = 0$, and if $\link \sigma$ is a homology sphere, $\chi (\link \sigma) - 1 = (-1)^{r - \dim \sigma-1}$.  Thus, by Equation (\ref{eq:5}),
	\begin{align*}
	\iota([M])
	&= \sum_{\sigma\not\in \partial M} (-1)^{\dim \sigma} \bigg((-1)^{r - \dim \sigma -1 }\bigg)\sigma 
		+\sum_{\sigma \in \partial M} (-1)^{\dim \sigma}\cdot 0\cdot \sigma\\
	&= \sum_{\sigma\not\in \partial M} (-1)^{r - 1}\sigma
	= (-1)^{r+1}\big([M] - [\partial M]\big).
	\end{align*}
\end{proof}

\begin{proof}[Proof of Lemma \ref{lem:generalcommutativity}]
By linearity, it suffices to prove that $\iota b_\Phi([\sigma]) = b_\Phi \iota([\sigma])$ for any simplex $\sigma$, where $[\sigma] = \sum_{\tau\subseteq\sigma}\tau$ is the sum of the simplices contained in the manifold $\sigma$. Since $\Phi$ is a subdivision method, $b_\Phi([\sigma])$ will also be the sum of the simplices contained in some manifold $\Phi(\sigma)$ of dimension $\dim\sigma$. Also, $b_\Phi([\partial\sigma]) = [\partial\Phi(\sigma)]$.

The result now follows from Lemma~\ref{lem:manifold formula}:
\begin{align*}
b_\Phi \iota([\sigma]) &= b_\Phi\big((-1)^{\dim \sigma+1}([\sigma] - [\partial\sigma])\big) \\
&= (-1)^{\dim \sigma+1}(b_\Phi([\sigma]) - b_\Phi([\partial\sigma])) \\
&= (-1)^{\dim \sigma+1}([\Phi(\sigma)] - [\partial\Phi(\sigma)]) \\
&= \iota([\Phi(\sigma)]) \\
&= \iota b_\Phi([\sigma]).
\end{align*}
\end{proof}

\begin{theorem}\label{thm:generalsym} %(see Theorem \ref{thm:generalsymmetry})
For any dimension $n$ and any subdivision method $\Phi$ which is nontrivial in dimension $n$, there exists a unique `limit polynomial' $p_{n,\Phi}(t)$, such that, for any $(d-1)$-dimensional simplicial complex $X$, the roots of $f^{\Phi^k(X)}(t)$ converge to the roots of $p_{n,\Phi}(t)$ as $k$ increases.  The roots of $p_{n,\Phi}(t)$ are invariant under the M\"{o}bius transformation $x \mapsto \frac{-x}{x+1}$.
\end{theorem}

\begin{proof}
The proof of this theorem is exactly the same as the proof of Theorem \ref{thm:generalsym} in the case of barycentric subdivision.  The key observations there were that $\iota b_\Phi = b_\Phi\iota$ and that there exists a unique eigenvector for the maximal eigenvalue of the matrix realizing the effect of subdivision on $f$-vectors.  That this eigenvector is unique in general follows from $\Phi$ being nontrivial in dimension $n$, and is left as an exercise for the reader.
\end{proof}

\begin{remark}
Since the above interpretation is on the level of formal sums of simplices, the most natural context in which to study it seems to be the Stanley-Reisner ring $\mathbb K[X]$, defined for any simplicial complex $X$ and any field $\mathbb K$. A good introduction to these rings can be found in \cite{StanleyCCA}, where some properties of the Stanley-Reisner ring of a subdivision of a simplicial complex are explored. This brings us to ask the following questions.
\end{remark}

\begin{question}
Is there a (multi-)complex in each dimension whose $f$-polynomial is related to the limit polynomials $\pol_\infty(t)$ or $\qol_\infty(t)$?
\end{question}

More generally: 

\begin{question}
Is there a geometric interpretation of the coefficients or the roots of $\pol_\infty(t)$ (equivalently, $\qol_\infty(t)$)?
\end{question}

%%%%%%%%%%%%%%%%%%%%%%%%%%%%%%%%%%%%%%%%%%%%%%%%%%%
\section{Computations}\label{sec:computations}

We finish this paper by computing explicit values for the limit roots up to $d = 10$.  As observed in Theorem \ref{thm:facts}, to compute $p_d$ we need to compute the matrix $P_d^{-1}$.  To do so, we first compute a more explicit expression for $\Lambda_d$.

Recall that $\Lambda_d = [\lambda_{i,j}]$.  

\begin{lemma}
	$$
	\lambda_{i-1,j}  = \sum_{l=0}^j (-1)^{j-l} { j \choose l } l^i.
	$$
\end{lemma}

\begin{proof}

Starting with Lemma \ref{lem:single}, it follows from Lemma \ref{lem:foij} that:

\begin{align*}
\lambda_{i-1,j} &= \sum_{k=0}^i (-1)^k{i \choose k } f_j^{\sigma'_{i-k}}\\
&= \sum_{k=0}^i (-1)^k{i \choose k }\Delta^j \{(1+l)^{i-k}\}_l\\
&= \sum_{k=0}^i (-1)^k{i \choose k }\sum_{l=0}^j(-1)^l { j \choose l } (1+l)^{i-k}\\
&= \sum_{l=0}^j (-1)^{j-l} { j \choose l } \sum_{k=0}^i (-1)^k{i \choose k } (1+l)^{i-k}\\
&= \sum_{l=0}^j (-1)^{j-l} { j \choose l } (1+l-1)^{i}\\
&= \sum_{l=0}^j (-1)^{j-l} { j \choose l } l^i.
\end{align*}

\end{proof}

Now that we have the coefficients of $\Lambda_d$, we may compute the coefficients of the diagonalizing matrix $P_d$ by computing the eigenvectors of $\Lambda_d$.  We may then compute $P_d^{-1}$ using standard inversion techniques.  Note that it is possible to compute the entries of both $P_d$ and $P_d^{-1}$ explicitly in terms of the entries of $\Lambda_d$ and hence in terms of sums and products of integers.  To give the reader an idea of the numerical consequences of these calculations, we report below the results obtained using a standard symbolic computation program.

For any $k \geq d$, the $d^{th}$ row of $P_k^{-1}$ does not depend on $k$, and gives the coefficients of $q_d(t)$.  Thus, we present here the matrix $P_{10}^{-1}$:
$$
\left[\begin{array}{ccccccccccc}
1&0&0&0&0&0&0&0&0&0&0\\[.5em]
0&1&0&0&0&0&0&0&0&0&0\\[.5em]
0&
1&1&0&0&0&0&0&0&0&0\\[.5em]
0& 
\frac{1}{2}& 
\frac{3}{2}& 1&0&0&0&0&0&0&0\\[.5em]
0& 
\frac{2}{11}& 
\frac{13}{11}& 
2&1&0&0&0&0&0&0\\[.5em]
0& 
\frac{1}{19}& 
\frac{25}{38}& 
\frac{40}{19}& 
\frac{5}{2}& 1&0&0&0&0&0\\[.5em]
0& 
\frac{132}{10411}& 
\frac{3004}{10411}& 
\frac{45}{29}& 
\frac{95}{29}& 
3&1&0&0&0&0\\[.5em]
0& 
\frac{90}{34399}& 
\frac{3626}{34399}& 
\frac{61607}{68798}& 
\frac{245}{82}& 
\frac{385}{82}& 
\frac{7}{2}& 1&0&0&0\\[.5em]
0& 
\frac{15984}{33846961}& 
\frac{12351860}{372316571}& 
\frac{7924}{18469}& 
\frac{39221}{18469}& 
\frac{56}{11}& 
\frac{70}{11}& 
4&1&0&0\\[.5em]
0& 
\frac{983304}{12980789207}& 
\frac{119432466}{12980789207}& 
\frac{2296176994}{12980789207}& 
\frac{536193}{429266}& 
\frac{919821}{214633}& 
\frac{567}{71}& 
\frac{588}{71}& 
\frac{9}{2}& 1&0\\[.5em]
0& 
\frac{1345248918720}{123031432784730871}& 
\frac{281136722386176}{123031432784730871}& 
\frac{4358731100}{67808366729}& 
\frac{42780833020}{67808366729}& 
\frac{1335075}{448471}& 
\frac{3478503}{448471}& 
\frac{1050}{89}& 
\frac{930}{89}& 
5&1
\end{array}
\right]
$$

The roots of $p_d(t)$ are, for $d \leq 10$, approximated by:
$$\left\{
\begin{array}{ccccccccccc}
d=2:&	-1\\
d=3:&	-2&	-1\\
d=4:&	-4.1861&	-1.3139&	-1\\
d=5:&	-8.3642&	-2	&	-1.1358&	-1\\
d=6:&	-16.096&	-1.4706&	-3.1252&	-1.0662&	-1\\
d=7:&	-30.121&	-4.8761&	-2	&	-1.2570&	-1.0343&	-1\\
d=8:&	-55.208&	-7.5398&	-2.7664&	-1.5661&	-1.1529&	-1.0185&	-1\\
d=9:&	-99.626&	-11.537&	-3.8404&	-2	&	-1.3521&	-1.0949&	-1.0101&	-1\\
d=10:&	-177.68&	-17.474&	-5.3206&	-2.5830&	-1.6317&	-1.2315&	-1.0607&	-1.0057&	-1
\end{array}
\right.$$

We see that these roots are symmetric about the point $-2$, with respect to the M\"{o}bius transformation $x \mapsto \frac{-x}{x+1}$.  In other words, if $a$ is a root of $p_d(t)$, then so is $\frac{-a}{a+1}$, where note $-2$ is fixed by this transformation. The symmetry is more apparent in the (linear) symmetry about $-\frac{1}{2}$ exhibited by the roots of $q_d(t)$, which are the reciprocals of the roots of $p_d(t)$.  The roots of $q_d(t)$ are, for $d \leq 10$, approximated by:
$$\left\{
\begin{array}{ccccccccccc}
d=2:&	-1&	 0\\
d=3:&	-1&	-.5		&	0\\
d=4:&	-1&	-.76112	&	-.23888	&	0\\
d=5:&	-1&	-.88044	&	-.5		&	-.11956	&	0\\
d=6:&	-1&	-.93787	&	-.68002	&	-.31998	&	-.06213	&	0\\
d=7:&	-1&	-.96680	&	-.79492	&	-.5		&	-.20508	&	-.03320	&	0\\
d=8:&	-1&	-.98189	&	-.86737	&	-.63852	&	-.36148	&	-.13263	&	-.01811&	0\\
d=9:&	-1&	-.98996	&	-.91332	&	-.73961	&	-.5		&	-.26039	&	-.08668&	-.01004&	0\\
d=10:&	-1&	-.99437	&	-.94277	&	-.81205	&	-.61285	&	-.38715	&	-.18795&	-.05723&	-.00563&	0
\end{array}
\right.$$

Further computations of the coefficients and roots have been carried out by H\"useyin \"Ozoguz at the University of Bremen.  Interesting patterns emerge.  For instance, \"Ozoguz has noted (without proof) that the denominators of the reduced fraction representations of the coefficients are all square-free, and that certain ratios between roots appear to converge as $d$ increases.  We leave the study of these phenomena to future work.

%%%%%%%%%%%%%%%%%%%%%%%%%%%%%%%%%%%%%%%%%%%%%%%%%%%

\bibliographystyle{plain}
\bibliography{DPS2}

\end{document}